\documentclass[11pt]{article}
\usepackage[a4paper,margin=1in]{geometry}
\usepackage{amsmath,amssymb,amsthm,mathtools}
\usepackage[hidelinks]{hyperref}
\hypersetup{
  pdftitle={A power-saving asymptotic for a Mobius-weighted GCD quadratic form},
  pdfauthor={Ahmadreza Azimifard}
}
\newtheorem{theorem}{Theorem}[section]
\newtheorem{lemma}[theorem]{Lemma}

\theoremstyle{remark}

\newcommand{\eps}{\varepsilon}
\newcommand{\1}{\mathbf 1}

\title{A power-saving asymptotic for a M\"obius-weighted GCD quadratic form}
\author{Ahmadreza Azimifard\thanks{Harmonic Research \& Technologies, LLC, New York, NY, USA.
  Email: \href{mailto:afard@harmonicrt.com}{afard@harmonicrt.com}}}
\date{}

\begin{document}
\maketitle

\begin{abstract}
For the M\"obius-weighted quadratic form
\[
 T(N)=\sum_{m,n\leq N}\frac{\mu(m)\mu(n)(m,n)^2}{mn},
\]
we prove $T(N)=cN+O(\sqrt N)$, where $c>0$ is given by an absolutely
convergent series multiplied by a convergent Euler product.
The Jordan-totient decomposition gives positivity, and a uniform mean
value for a squarefree weight gives the error after interchanging finite
sums.
\end{abstract}

\section{Statement of the result}

The matrices with entries $(m,n)^2/(mn)$ form a classical family of positive
semidefinite GCD matrices.  Here we evaluate their quadratic forms on the
M\"obius vector.  More precisely, for a positive integer $N$, let
\begin{equation}\label{eq:T}
 T(N)=\sum_{m,n\leq N}\frac{\mu(m)\mu(n)(m,n)^2}{mn},
\end{equation}
where $\mu$ denotes the M\"obius function and $(m,n)$ the greatest common divisor
of $m$ and $n$.  We shall identify the linear term in $T(N)$ and obtain a
power saving in the remainder.

The constant is most naturally expressed using
\begin{equation}\label{eq:C0R}
 C_0=\prod_p(1-p^{-1})(1+p^{-1}-p^{-3}),
 \qquad
 R(r)=\prod_{p\mid r}\frac{p^3}{p^3+p^2-1},
\end{equation}
and
\begin{equation}\label{eq:c}
 c=C_0\sum_{q\geq1}\frac{A_q}{q(q+1)},\qquad
 A_q=\sum_{k,l\leq q}\frac{\mu(k)\mu(l)R(kl)}{kl}.
\end{equation}
Each factor in the product for $C_0$ is positive and equals $1+O(p^{-2})$,
so the product converges to a positive number.  Also, $0<R(r)\leq1$.

\begin{theorem}\label{thm:main}
The series defining $c$ in \eqref{eq:c} converges absolutely, $c>0$, and
$A_q\geq0$ for every positive integer $q$. Moreover,
\[
 T(N)=cN+O(\sqrt N).
\]
\end{theorem}

The proof uses the square decomposition to determine the sign of the
constant. Expanding the squares and applying the uniform mean value
before grouping the pairs gives the square-root remainder.

Throughout, $J_2(n)=n^2\prod_{p\mid n}(1-p^{-2})$ denotes the second Jordan
totient, $\tau(n)$ the number of positive divisors of $n$, and
$H_q=\sum_{k=1}^q k^{-1}$.  Sums over arithmetic variables run over positive
integers unless stated otherwise.  Implied constants are absolute except for
the dependence indicated by a subscript.

\section{An exact positive decomposition}

We use the standard identity $n^2=\sum_{d\mid n}J_2(d)$ in the form
\begin{equation}\label{eq:jordan-kernel}
 \frac{(m,n)^2}{mn}=\sum_{d\mid(m,n)}\frac{J_2(d)}{mn}.
\end{equation}
This turns the quadratic form into a sum of nonnegative terms.

\begin{lemma}\label{lem:square}
For every positive integer $N$,
\begin{equation}\label{eq:square}
 T(N)=\sum_{d\leq N}\frac{J_2(d)}{d^2}
       \left(\sum_{k\leq N/d}\frac{\mu(dk)}k\right)^2.
\end{equation}
In particular, $T(N)\geq0$.
\end{lemma}

\begin{proof}
Substitute \eqref{eq:jordan-kernel} into \eqref{eq:T} and sum first over
the common divisor $d$.  Writing $m=dk$ and $n=dl$ leaves two identical
inner sums, whose product is the square in \eqref{eq:square}.
\end{proof}

The M\"obius factors satisfy
\begin{equation}\label{eq:mu-factor}
 \mu(dk)\mu(dl)=\mu(d)^2\mu(k)\mu(l)\1_{(d,kl)=1}.
\end{equation}
Indeed, $\mu(ab)=\mu(a)\mu(b)\1_{(a,b)=1}$ for positive integers $a,b$.
Apply this identity to $dk$ and $dl$, and combine the two coprimality conditions.

\section{A uniform mean value}

The coprimality condition in \eqref{eq:mu-factor} leads us to the functions
\begin{equation}\label{eq:fr}
 f_r(d)=\mu(d)^2\prod_{p\mid d}(1-p^{-2})\1_{(d,r)=1}.
\end{equation}
We use a squarefree convolution, as in the proof of the coprimality
estimate of Robert and Tenenbaum \cite[Sec.~10, Eq.~(10.1)]{RT2013}.

\begin{lemma}[uniform mean value]\label{lem:mean}
Uniformly for integers $r\geq1$ and real $X\geq1$,
\begin{equation}\label{eq:mean}
 \sum_{d\leq X}f_r(d)=C_0R(r)X+O\bigl(\tau(r)^2\sqrt X\bigr).
\end{equation}
\end{lemma}

\begin{proof}
Define the multiplicative function $h_r$ by
\[
 h_r(p^j)=
 \begin{cases}
 (-1)^j p^{-2},&p\nmid r,\\
 (-1)^j,&p\mid r,
 \end{cases}
 \qquad j\geq1.
\]
Multiplying the local series for $h_r$ by $1+z$ gives
$1+f_r(p)z$, so $f_r=\mu^2*h_r$, with $*$ denoting Dirichlet convolution.
At $z=p^{-1/2}$ the product of the absolute local series satisfies
\[
 \sum_{e\geq1}\frac{|h_r(e)|}{\sqrt e}
 \leq \prod_p\left(1+\frac{p^{-5/2}}{1-p^{-1/2}}\right)
        \prod_{p\mid r}\frac{1}{1-p^{-1/2}}
 \ll 4^{\omega(r)}\leq\tau(r)^2,
\]
where $\omega(r)$ counts distinct prime factors and the first product
converges.

The identity $\mu^2(d)=\sum_{b^2\mid d}\mu(b)$ gives
$\sum_{d\leq Y}\mu^2(d)=Y/\zeta(2)+O(\sqrt Y)$ for $Y\geq1$.
Summing the convolution up to $X$ yields the main coefficient
$\zeta(2)^{-1}\sum_{e\geq1}h_r(e)/e$. Indeed, the accumulated error and
the main-series tail are bounded by a constant times
$\sqrt X\sum_e |h_r(e)|/\sqrt e$.
The coefficient has Euler factors
\[
 (1-p^{-2})\frac{1+f_r(p)/p}{1+1/p}
 =(1-1/p)(1+f_r(p)/p).
\]
Their product is $C_0R(r)$, as claimed.
\end{proof}

\begin{lemma}\label{lem:positive-constant}
For every $q\geq1$, we have $0\leq A_q\leq H_q^2$.  The series in
\eqref{eq:c} converges absolutely, and $c\geq C_0/2>0$.
\end{lemma}

\begin{proof}
Fix $q$. Expanding a square and using
\eqref{eq:mu-factor} and Lemma~\ref{lem:mean}, we obtain
\begin{align*}
 S_q(X)&:=\sum_{d\leq X}\frac{J_2(d)}{d^2}
              \left(\sum_{k\leq q}\frac{\mu(dk)}k\right)^2\\
 &=\sum_{k,l\leq q}\frac{\mu(k)\mu(l)}{kl}
              \sum_{d\leq X}f_{kl}(d)
   =C_0A_qX+O_q(\sqrt X).
\end{align*}
Since $S_q(X)\geq0$, division by $X$ followed by $X\to\infty$ gives
$C_0A_q\geq0$.  Thus $A_q\geq0$ because $C_0>0$.  The triangle inequality
and $R(kl)\leq1$ also give $A_q\leq H_q^2$.  Absolute convergence follows
from $H_q\leq1+\log q$ and convergence of
$\sum_{q\geq1}(1+\log q)^2/q^2$.  Finally, $A_1=R(1)=1$, so the first
term of the nonnegative series gives $c\geq C_0/2$.
\end{proof}

\section{The asymptotic estimate}

\begin{proof}[Proof of Theorem~\ref{thm:main}]
Lemma~\ref{lem:positive-constant} proves convergence and positivity.
Expand \eqref{eq:square}, use \eqref{eq:mu-factor}, and sum in the opposite
order. With $M=\max(k,l)$ this gives
\[
 T(N)=\sum_{k,l\leq N}\frac{\mu(k)\mu(l)}{kl}
                \sum_{d\leq N/M} f_{kl}(d).
\]
By Lemma~\ref{lem:mean},
\[
 T(N)=NC_0\sum_{k,l\leq N}\frac{\mu(k)\mu(l)R(kl)}{klM}
       +O\left(\sqrt N\sum_{k,l\leq N}
                    \frac{\tau(kl)^2}{kl\sqrt M}\right).
\]
For any fixed $0<\eta<1/4$, the standard divisor bound gives
$\tau(kl)^2\ll_\eta(kl)^\eta$. Grouping pairs by $M$ bounds the error
series by a constant times
$\sum_{M\geq1}M^{-3/2+\eta}\sum_{l\leq M}l^{-1+\eta}$,
which converges because the inner sum is $O_\eta(M^\eta)$.
Since $R(kl)\leq1$, extending the main series to all pairs costs at most
\[
 2NC_0\sum_{M>N}\frac{H_M}{M^2}=O(\log(2N)).
\]
Finally, the identity
$\sum_{q\geq M}1/(q(q+1))=1/M$ gives the absolutely convergent expression
\[
 c=C_0\sum_{k,l\geq1}\frac{\mu(k)\mu(l)R(kl)}{kl\max(k,l)}.
\]
This proves the square-root estimate.
\end{proof}

\section{Formal verification and numerical checks}

The ancillary Lean project verifies convergence and positivity with
an error estimate of the form
\[
 |T(N)-Nc|\leq C_\eps N^{2/3+\eps}\qquad(N\geq8,\ \eps>0),
\]
as in the earlier version. The formal constant uses a convergent
coefficient series. Its identification with the displayed Euler product
is not a separate exported theorem. The square-root estimate in the
present version is proved above on paper.
The formal project used Lean~4.28.0 and Mathlib~4.28.0. The checked
headline declaration uses only \texttt{propext},
\texttt{Classical.choice}, and \texttt{Quot.sound}. The ancillary README
gives the correspondence and build instructions.

For a numerical comparison, floating-point evaluation of the square
decomposition gives
\[
 \frac{T(10^6)}{10^6}\approx0.3311235851.
\]
Truncating the series in \eqref{eq:c} at $q=5000$ and the product for $C_0$
at primes $p\leq2\cdot10^6$ gives approximately $0.3311220400$.
These computations are consistent with the theorem, but no rigorous error
enclosure is supplied for either decimal.  They play no part in the proof.

\section{Relation to prior work}

The GCD kernel used here and its interpretation as an $L^2$ Gram kernel are
classical; see Lindqvist and Seip~\cite{LS1998}.  Aistleitner, Berkes, Seip and
Weber~\cite{ABSW2015} study GCD matrices and their spectral norms in connection
with series of dilated functions.  Hilberdink~\cite{Hilberdink2014} develops a
related framework on the group of squarefree integers, and Ehm~\cite{Ehm2024}
studies related Gram matrices and series.  A Jordan-weighted square identity
also appears in Montgomery and Vaughan's treatment of Farey discrepancy;
see Appendix~F of~\cite{MVII}.

For a positive integer $X$, let $E_X$ denote averaging over $[0,X]$.
Chavez~\cite{Chavez2024} obtains $-9/(2\pi^2)+O(1/X)$ for the
off-diagonal sum of $\mu(m)\mu(n)E_X(\{mx\}\{nx\})/(mn)$ over distinct
positive $m,n$. Its coefficient weights, uncentred factors and omitted
diagonal distinguish this quantity from $T(N)$.

Caceres~\cite[Thm.~10.2 and Rem.~10.3]{Caceres2026} studies exactly
$V(D)=T(D)$, decomposes it into nonnegative terms and identifies $V(D)/12$
with the long-run mean square of $\sum_{d\le D}\mu(d)B_1(x/d)$,
where $B_1(t)=\{t\}-1/2$ away from integers.
He reports $V(D)/D\approx0.33114$ at $D=2\cdot10^6$ and leaves the
existence of its limit open. Theorem~\ref{thm:main} gives that limit $c$
with an error term. Our numerical value $T(10^6)/10^6\approx0.3311236$
is consistent with his observation.

Dress, Iwaniec and Tenenbaum~\cite{DIT1983} prove that
$S(z)=\sum_{m,n\le z}\mu(m)\mu(n)/[m,n]$ tends to a positive limit,
where $[m,n]$ is the least common multiple.
Since $1/[m,n]=(m,n)/(mn)$ and $(m,n)=\sum_{d\mid(m,n)}\varphi(d)$,
Euler's totient $\varphi$ and the second Jordan totient $J_2$ give
\[
 S(z)=\sum_{d\le z}\varphi(d)
       \left(\sum_{\substack{m\le z\\d\mid m}}\frac{\mu(m)}m\right)^2,
 \qquad
 T(N)=\sum_{d\le N}J_2(d)
       \left(\sum_{\substack{m\le N\\d\mid m}}\frac{\mu(m)}m\right)^2 .
\]
Their proof combines square decomposition, coprimality-dependent mean
values and quotient blocks. The present form uses the second Jordan
totient and has linear growth. The squarefree coprimality estimate of
Robert and Tenenbaum \cite[Eq.~(10.1)]{RT2013} also underlies later work of
de la Bret\`eche, Dress and Tenenbaum \cite[Eq.~(2.1)]{BDT2020} on the
LCM-weighted form.

For $\psi(t)=\{t\}-1/2$, Franel's identity in the Farey literature
\cite{Franel1924,Landau1924} reads
$\int_0^1\psi(mx)\psi(nx)\,dx=(m,n)^2/(12mn)$.
Thus Theorem~\ref{thm:main} gives a mean-square asymptotic through
\[
 T(N)=12\int_0^1\left(\sum_{n\le N}\mu(n)\psi(nx)\right)^2dx .
\]

\section*{Use of computational tools}
Automated research tools contributed to the initial derivation and draft.
Generative AI assisted with literature comparisons, proof development and
revision.

\bibliographystyle{plain}
\bibliography{refs}

\end{document}